\documentclass[12pt]{article}
\usepackage{graphicx}

\font\tendb=msbm10 at 12pt \font\sevendb=msbm10 at 9pt
\font\fivedb=msbm10 at 7pt
\newfam\dbfam
\textfont\dbfam=\tendb \scriptfont\dbfam=\sevendb
\scriptscriptfont\dbfam=\fivedb
\def\db{\fam\dbfam\tendb}
\font\eufm=eufm10\font\eufms=eufm10\font\eufmss=eufm10\newfam\eufam
\textfont\eufam=\eufm\scriptfont\eufam=\eufms\scriptscriptfont\eufam=\eufmss

\font\tendbb=msbm10 at 12pt \font\sevendbb=msbm7 at 9pt
\font\fivedbb=msbm5 at 6pt
\newfam\dbbfam
\textfont\dbbfam=\tendbb \scriptfont\dbbfam=\sevendbb
\scriptscriptfont\dbbfam=\fivedbb

 \def \Z {{\db Z}}
 \def \R {\hbox{\db R}}
 
 \def \C {\hbox{\db C}}
 
 \def \S {S^{3}}

\def \fin {\hfill \framebox(7,7) }
 
\font\tenMmm=eusm10 at 12pt
\newfam\Mmmfam
\textfont\Mmmfam=\tenMmm

\def\illu #1 by #2 (#3){
  \vbox to #2{
    \hrule width #1 height 0pt depth 0pt
    \vfill
    \special{illustration #3} 
    }
  }

\begin{document}

\null \vspace{1cm}

\begin{center}  \textbf{\large {Strong periodicity of links
and the coefficients of the Conway polynomial}}\\
Nafaa Chbili\footnote{ Supported by a
fellowship from the COE
 program "Constitution of wide-angle mathematical basis focused on
 knots",
  Osaka City University. The author would like to express his thanks and
  gratitude to
  Akio Kawauchi for his kind hospitality.}\\
\begin{footnotesize}
Osaka City University Advanced Mathematical Institute \\
Sugimoto 3-3-138, Sumiyoshi-ku 558 8585 Osaka, Japan. \\
 E-mail: chbili@sci.osaka-cu.ac.jp\\
  \end{footnotesize}
\end{center}

 \begin{footnotesize}
               \begin{abstract} Przytycki and Sokolov proved that
 a three-manifold admits a semi-free action of the finite cyclic group of
 order $p$ with a circle
 as the set of fixed points if and only if $M$ is obtained from the
 three-sphere by surgery along a strongly $p-$periodic link $L$. Moreover,
 if the quotient three-manifold is an integral homology sphere, then
 we may assume that $L$ is orbitally separated. This paper studies the
behavior  of the coefficients of the Conway polynomial of such a
link. Namely, we prove that if $L$ is a strongly $p$-periodic
orbitally separated link and  $p$ is an odd prime, then
 the coefficient $a_{2i}(L)$ is congruent to zero modulo $p$
 for all $i$ such that $2i<p-1$.    
 \end{abstract}
               \end{footnotesize}

\begin{center} \textbf{1- Introduction}
\end{center}
 Strongly periodic links have been introduced by Przytycki and
Sokolov \cite{PS} as a class of links providing equivariant
surgery presentations of cyclic covers branched along knots
\emph{i.e.} \emph{periodic three-manifolds}. This equivariant
surgery presentation plays a key role in the study of the
Casson-Walker-Lescop invariant as well as the quantum invariants
of periodic three-manifolds \cite{Ch1, Ch2, Ch3, CL, BGP}. Let $p
\geq 2$ be an integer, a link $L$ in $\S$ is said to be
\emph{$p$-periodic} if there exists a periodic transformation $h$
of order $p$ such that $fix(h)\cong S^1, h(L)=L$ and $fix(h)\cap
L=\emptyset$. By the positive solution of the Smith conjecture
\cite{BM} we may assume that $h$ is a rotation by a $2\pi/p$ angle
around the $z$-axis. A $p$-periodic link is said to be
\emph{strongly $p$-periodic} if and only if every component of the
quotient link $\overline L$ has linking number zero modulo $p$
with the axis of the rotation. The link $L$ is said to be
\emph{orbitally separated} (we write OS for short) if the quotient
link is algebraically split. Przytycki and Sokolov \cite{PS}
proved that a three-manifold is periodic with prime period $p$ if
and only if $M$ is obtained from $\S$ by surgery along a strongly
$p-$periodic link $L$. Moreover, It can be easily seen that if the
quotient manifold $\overline M$ is an integral homology sphere
then we may assume that $L$ is orbitally separated.\\
 Murasugi \cite{Mu} established  a congruence formula
for the Alexander polynomial of   periodic knots. Later,
 different approaches   have been  used to give alternative  proofs to Murasugi's
 result \cite{Bu,DL,Hi}. In \cite{Mi} and \cite{Sa}, similar
 results have been  proved for the multi-variable Alexander
 polynomial of periodic links.\\
 The Conway polynomial is a renormalization  of the
 Alexander polynomial. If $L$ is a link with $n$
 components then the Conway polynomial of $L$ is of the form $
\nabla_{L}(z)=z^{n-1}(a_0(L) + a_2(L) z^2+ \dots + a_{2m}(L)
z^{2m})$, where coefficients $a_i(L)$ are integers. We know that
the coefficient $a_0(L)$ depends only on the linking matrix of
$L$. The second coefficient $a_2(L)$ is related to the
Casson-Walker-Lescop invariant by the global surgery formula
introduced by Lescop \cite{Les}. The study of the  coefficients of
the Conway polynomial is of interest  partially because they
 are Vassiliev invariants of suitable orders.\\
Let $p$ be an odd prime and let  $L$ be an OS strongly
$p$-periodic link. In \cite{Ch3} we proved that the coefficient
$a_2(L)$ is congruent to zero modulo $p$. Using this fact and
Lescop's surgery formula  we have been able to prove that the
 Casson-Walker-Lescop invariant of a periodic three-manifold is equal, modulo $p$, to
 the product of the signature of the surgery presentation by the order of the first
 homology group over 8.\\
The purpose of this paper is to extend the result proved in
\cite{Ch3} for the coefficient $a_2$ to other coefficients of the
Conway
polynomial. Here is our Main Theorem\\

 \textbf{Main Theorem.} {\sl Let $p$ be an odd  prime and let  $L$ be an OS strongly
$p$-periodic link, then $a_{2i}(L)$ is congruent to zero modulo
$p$ for all $i$ such that $2i<p-1$.}\\

The Main Theorem does not hold for $p=2$. A counterexample is
given by considering the Hopf link $H$. Obviously, $H$ is an OS
strongly 2-periodic link. The Conway polynomial of $H$ is
$\nabla_{H}(z)=z$. This means that the coefficient $a_{0}(H)=1$,
hence it is  not null modulo 2.\\

 Here is an outline of the paper. Section 2  introduces strongly
periodic links and discusses the notion of  equivariant crossing
changes. Section 3 reviews some properties of the Conway
polynomial needed in the sequel.
Finally, Section 4 is devoted to the proof of the Main Theorem.\\

\begin{center} \textbf{{2- Strongly periodic links
and equivariant crossing changes}}
\end{center}

In this section, we first  define strongly periodic links. Then we
introduce some of their properties. We close the section with a
list discussing the different types of  equivariant crossing changes.\\

{\bf Definition 2.1.} {\sl { Let $p\geq 2$ be an integer. A link
$L$ of $\S$ is said to be $p$-periodic if and only if there exists
an orientation-preserving  auto-diffeomorphism $h$ of  $\S$
such that:\\
\begin{tabular}{rl}
&1- Fix($h$) is homeomorphic to the  circle $S^{1}$,\\
&2- the link  L is disjoint from  Fix($h$),\\
&3- $h$ is of order $p$,\\
&4- $h(L)=L$.
\end{tabular}

If $L$ is  periodic we will denote the quotient link by $\overline
L$.\\}}

The standard example of such a diffeomorphism is given as follows.
Let us consider  $\S$ as the  sub-manifold  of $\C ^2$ defined by
$\S=\{(z_{1},z_{2}) \in \C ^{2}; |z_{1}|^{2}+|z_{2}|^{2}=1\}$ and
$\varphi{_p}$ the following diffeomorphism:
$$
\begin{array}{cccl}
\varphi_{p}:& S^{3} & \longrightarrow & S^{3} \\
   & (z_{1},z_{2}) & \longmapsto & (e^{\frac{2i\pi}{p}}z_{1},z_{2}).
\end{array} $$
The set of fixed points for  $\varphi_{p}$ is the circle $\Delta =
\{(0,z_{2}) \in \C ^{2}; |z_{2}|^{2}=1\}$. If we identify $\S$
with $\R ^{3}\cup {\infty}$, $\Delta$ may be seen as the standard
$z$-axis.\\
Recall here that if the quotient link $\overline L$ is a knot then
the link $L$ may have more than one component. In general, the
number of components of $L$ depends on the linking numbers of the
components of $\overline L$ with the axis
of the rotation.\\
Recently, Przytycki and Sokolov \cite{PS} introduced the notion of
strongly periodic links as follows.\\

{\bf Definition 2.2.} {\sl  Let $p\geq 2$ be an integer. A
$p$-periodic link $L$ is said to be strongly $p$-periodic if and
only if one of the following conditions holds:\\
\begin{tabular}{rl}
&(i) The linking number of each component of $\overline L$ with
the axis $\Delta$ is congruent
 to 0  modulo $p$.\\
 &(ii) The group $\Z/{p\Z}$ acts freely on the set of components of
 $L$.\\
 &(iii) The number of components of $L$ is $p$ times greater than
 the number of components $
\overline L$.
\end{tabular}
} \\

 {\bf Remark 2.1.} According to condition (iii)  in the previous definition,
 a $p$-strongly
 periodic link $L$ has $p\alpha$ components, where $\alpha$ is the
 number of components of the quotient link $\overline L$. These
 $p\alpha$ components are divided into $\alpha$ orbits with
 respect to the free action of  $\Z/{p\Z}$ (condition (i)).\\
 Assume that  $\overline L=l_1\cup \dots \cup\l_{\alpha}$, there is a natural
cyclic order on each orbit of components of $L$. Namely,
$$L=l_1^1 \cup\dots \cup l_1^p \cup l_2^1 \cup\dots \cup l_2^p \cup \dots \cup l_{\alpha}^1 \cup\dots
\cup l_{\alpha}^p$$ where $ \varphi_p(l_i^t)=l_{i}^{t+1} \;
\forall 1\leq t\leq p-1$ and $\varphi_p(l_i^p)=l_{i}^{1}$, for all
$1\leq i \leq \alpha$.\\

 {\bf Definition 2.3.}  {\sl A link  $L$  in the three-sphere is said to be
  algebraically split if and only if the linking
number of any two components of $L$ is null.}\\

 {\bf Definition 2.4.} {\sl  Let $p\geq 2$ be an integer. A strongly $p$-periodic link is said to be
 orbitally separated if and only if the quotient link is
 algebraically split.}\\

 {\bf Remark 2.2.} A three-manifold $M$ is said to be $p-$periodic
 if the finite cyclic group of order $p$ acts semi-freely on $M$
 with a circle as the set of fixed points. By \cite{PS} and \cite{Sa1}, we know
  that a three-manifold is $p-$periodic if and only if $M$ is
  obtained from $\S$ by surgery along a
strongly $p-$periodic link.  If the quotient manifold is an
integral homology sphere, then we may choose the link $L$
orbitally separated. This is a consequence of the fact that the
quotient
manifold may be obtained by surgery along an algebraically split link.\\

 {\bf Remark 2.3.} For short, we shall use the term  OS link to refer to
  an orbitally separated link.
  If $L$ is a  strongly $p$-periodic OS link, then for all $s$ and $i\neq j$ we have
$\displaystyle\sum_{t=1}^p lk(l_i^s,l_j^t)=0$.\\

  \textbf{Definition 2.5.} {\sl Let $m$ be a positive integer. A $p-$periodic link
  $L$ is said to
  be of
type $m$ if $L=K_1 \cup K_1'\cup \dots \cup K_m\cup K_m'
\cup L'$   such that\\
\begin{tabular}{rl}
& (i) $K_i$ and $K_i'$  are invariant by the rotation for all $i$,
and $L'$ is an OS strongly
$p-$periodic link,\\
&(ii) Lk($K_i, \Delta) \equiv -$Lk$(K_i',\Delta)$ modulo $p$,
\\
 &(iii)
Lk$(\overline K_i,l)=-$Lk$(\overline {K_i'},l)$ for all components
$l$ of the quotient link $\overline {L'}$.
\end{tabular}

 By convention a $p$-periodic link  of type $0$ is a strongly
$p-$periodic OS link.}\\

Let $L$ be a $p$-periodic link in the three-sphere. Let $\overline
L$ be the factor link, so here we have $L=\pi^{-1}(\overline L)$,
where $\pi$ is the canonical surjection corresponding to the
action of the rotation on the three-sphere. Let $\overline L_+$,
$\overline L_-$ and $\overline L_0$ denote the three links which
are identical to $\overline L$ except near one crossing where they
are like in Figure 1. Now, let $L_{p+}:=\pi^{-1}(\overline L_+)$,
$L_{p-}:=\pi^{-1}(\overline L_-)$ and $L_{p0}:=\pi^{-1}(\overline
L_0)$. We define an {\em equivariant crossing change} as a change
 from $L_{p+}$ to $L_{p-}$ or vice-versa. If this crossing change
 involves
 two different components of the quotient link then we call it a
mixed equivariant crossing change. Otherwise, it is called a self
equivariant crossing change.\\

\textbf{Lemma 2.1.} {\sl Let $L=K_1 \cup K_1'\cup \dots \cup
K_m\cup K_m'
\cup L'$  be a  $p-$periodic link of type $m$.\\
1- If we change a crossing between $\overline K_i$ and  $\overline
{K_i'}$, then $L_{p0}$ has  $(\sharp L+p-2)$ components.\\
2- If we change a crossing between $\overline K_i$ or  $\overline
{K'_i}$ and $\overline K_j$ for $i \neq j$ then $L_{p0}$ has
either  $(\sharp L+p-2)$  or $(\sharp L-1)$ components.\\
3-  If we change a crossing between $\overline K_i$ or  $\overline
{K'_i}$ and a component of $\overline {L'}$, then $L_{p0}$ is of
type $m$ having $(\sharp L-p)$ components.\\
 4- If we change a crossing between two components of
$\overline {L'}$, then $L_{p0}$ is of type $m$ having $(\sharp L-p)$ components.\\
5- If we change a self-crossing in $\overline {L'}$,  then
$L_{p0}$ is either of type $(m+1)$ having $(\sharp L -p+2)$
components or a periodic link having  $(\sharp L +p)$ components.}\\
\emph{Proof.} The proof is straightforward by analyzing the
linking numbers of the components of $\overline L_{p0}$ with the
axis of the rotation.
 \fin \\

\begin{center}
{ \textbf{3- The Conway polynomial}}\\
\end{center}
The Conway polynomial $\nabla$ is an invariant of ambient isotopy
of oriented links which can be defined uniquely by the following:
$$\begin{array}{ll}
 {(i)}&\nabla_ {\bigcirc}(z)=1\\
{(ii)}&\nabla_{L_{+}}(z)-\nabla_{L_{-}}(z) =z\nabla_{L_{0}}(z),
\end{array}$$
where  $\bigcirc $ is the trivial knot, $L_{+}$,  $L_{-}$ and
$L_{0}$ are three oriented links which are identical except
 near one crossing where they look
like in the following figure:
\begin{center}
\includegraphics[width=8cm,height=2cm]{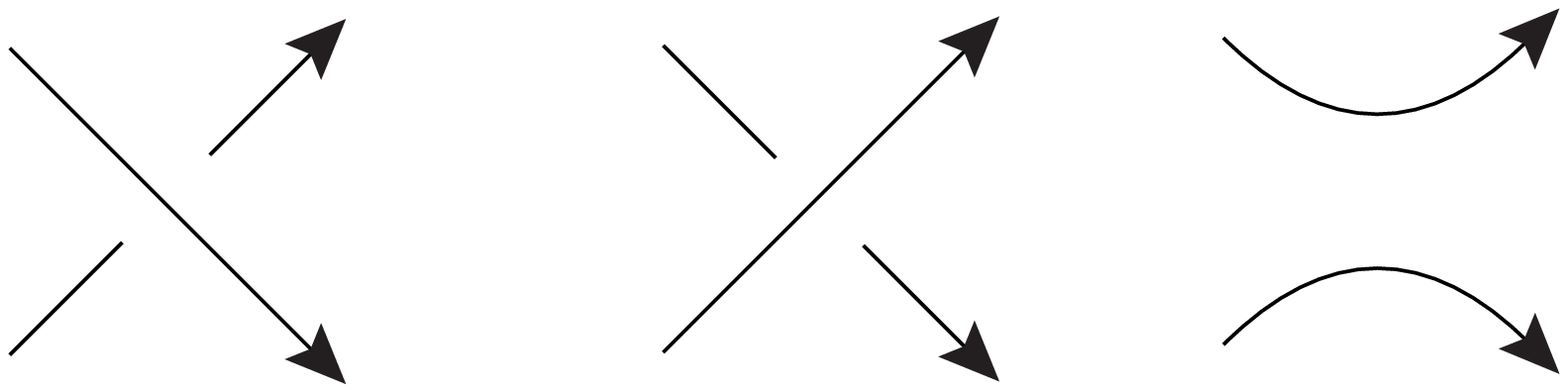}
\end{center}
$\hspace{150pt} L_{+}\hspace{80pt} L_{-}\hspace{70pt} L_{0}$
\begin{center} {\sc Figure 1.}
\end{center}

It is well known that if $L$ is a link with $n$
 components then the Conway polynomial of $L$ is of the form $
\nabla_{L}(z)=z^{n-1}(a_0(L) + a_2(L) z^2+ \dots + a_{2m}(L)
z^{2m})$,
where coefficients $a_i(L)$ are integers.\\
  Now, we shall
recall some properties of the Conway polynomial needed in the
sequel.\\
Let $L=l_1\cup \dots \cup l_n$ be an oriented $n$ component link
in the three-sphere. Let $l_{ij}=\mbox{Lk}(l_i,l_j)$ if $i \neq j$
and define $l_{ii}=-\displaystyle\sum_{j=1,j \neq i}^{n} l_{ij}$.
Define the linking matrix  ${\cal L}$ as ${\cal L}=(l_{ij})$. The
matrix $\cal L$ is symmetric and  each row adds to zero. Hence,
every cofactor of this matrix is the same. We refer the reader to
\cite{Ha}, \cite{Hoso} and \cite{Ho}, for
details about  the linking matrix.\\

 {\bf Theorem 3.1 \cite{Ho}.} {\sl Let $L$ an oriented link with
$n$ components. Then $a_0(L)={\cal L}_{ij}$, where  ${\cal
L}_{ij}$
is any cofactor of the linking matrix ${\cal L}$.}\\

It follows from this theorem that if $L$ is algebraically split
then the first coefficient of the Conway polynomial is zero.
Levine \cite{Le} proved that the coefficient of $z^i$ is null, for
all $i\leq 2n-3$, and
gave an explicit formula for the coefficient of $z^{2n-2}$.\\

 {\bf Proposition 3.2 \cite{Le}.} {\sl If  $L$ is algebraically split
with $n$ components, then $\nabla_{L}(z)$ is divisible by
$z^{2n-2}$.}\\


\begin{center}
\textbf{{4- Proof of the Main Theorem}} \end{center}

The idea of the proof is to use equivariant crossing changes to
reduce strongly periodic OS links to algebraically split links
without changing the first coefficients of the Conway polynomial
modulo $p$. We first introduce the following lemma which plays a
 crucial role in the proof.\\

\textbf{Lemma 4.1.} {\sl Let $p$ be a prime, then we have the
following congruence modulo $p$: }
$$
\nabla_{L_{p+}}(z)-\nabla_{L_{p-}}(z)\equiv z^p
\nabla_{L_{p0}}(z).
$$
\emph{Proof.} See \cite{Ch3}. \fin\\

\textbf{Remark 4.1.} Let $k$ be the number of components of
$L_{p+}$. Obviously, the link $L_{p-}$ has  $k$ components as
well. Let $s$ be the number of components of $L_{p0}$. According
to Lemma 4.1, we have the following congruence modulo $p$
$$z^{k-1}(a_0(L_{p+})+a_2(L_{p+})z^2+\dots)-z^{k-1}(a_0(L_{p-})+a_2(L_{p-})z^2+\dots)
\equiv z^{p+s-1}(a_0(L_{p0})+a_2(L_{p0})z^2+\dots).$$
Consequently, if $s \geq k-1$ then $a_{2i}(L_{p+}) \equiv
a_{2i}(L_{p-})$ modulo $p$ for all $i$ such that  $2i < p-1.$\\

\textbf{Lemma 4.2.} {\sl Let $L$ be a $p-$periodic link of type
$m$. If $L$ is algebraically split,
then $a_{2i}(L)$ is congruent to zero modulo $p$ for all $i$ such that $2i < p-1$.}\\

\emph{Proof.} We know that $L$ has $2m+p\alpha$ components. If
$\alpha>0$ then the number of components of $L$ is greater than or
equal to $p$. Using Levine's theorem we should be able to
conclude that $a_{2i}$ is
zero modulo $p$ for all $i$ such that  $2i < p-1$.\\
If $\alpha=0$ then $L=K_1 \cup K_1'\cup \dots \cup K_m\cup K_m'$.
The idea here is to make equivariant crossing changes until
getting a split link, without changing the coefficient $a_{2i}$
modulo $p$. We shall focus on equivariant  crossing changes
involving $K_1$ and the other components.  If we change a crossing
between $\overline K_1$ and $\overline K_1'$ then the link
$L_{p0}$ has $2(m-1)+p$ components. However, if we change a
crossing between $\overline K_1$  and $\overline K_j$ (resp.
$\overline K_j'$) for $1 \neq j $ then the number of components of
$L_{p0}$ is either $2(m-1)+p$ or $2(m-1)+1$. A simple computation
shows that in  all cases Remark 4.1 implies that $a_{2i}(L_{p+})$
is congruent to $a_{2i}(L_{p-})$ modulo $p$ for $2i<p-1$. It can
be easily seen that using these crossing changes we can put  $K_1$
over the rest of components of $L$ to get a split link. Since $L$
has the same coefficient $a_{2i}$ modulo $p$ than a split link,
then $a_{2i}(L)$ is null modulo $p$ for all $i$ such that $2i <
p-1$.
This ends the proof.\fin \\
\\

\textbf{Lemma 4.3.} {\sl If $L$ is of type $m>0$, then
$a_{0}(L)$ is null modulo $p$.}\\

\emph{Proof.} The coefficient $a_0(L) $ is computed using the
Hoste formula. Recall that $L$ can be written as follows:
$$L=K_1
\cup K'_1\cup \dots \cup K_m\cup K'_m  \cup l_{1}^1 \cup\dots \cup
l_{1}^p \cup \dots \cup l_{\alpha}^1 \cup\dots \cup
l_{\alpha}^p.$$
 It is easy
to see that the linking number of any couple of two invariant
components of $L $ is null modulo $p$. Moreover for such a
component $K$ we have:
$$
Lk(K,l_j^1)=Lk(K,l_j^2)=\dots =Lk(K,l_j^p)\;\; \forall 1\leq j
\leq \alpha.$$
 By Definition 2.5, (iii) we have Lk$(\overline K_i,l)=-$Lk$(\overline {K'_i},l)$ for
 all components $l$ of the quotient link $\overline {L'}$, where
  $L'$ is the strongly periodic part
 of $L$. Thus,
 $\sum_{t=1}^{p} \mbox{Lk}(K_i,l_j^t))=-\sum_{t=1}^{p} \mbox{Lk}(K_i',l_j^t))$
  for  all $1\leq
i \leq m$ and $1 \leq j \leq \alpha$. Notice that in each of  the
two summations above  the $p$ terms are equal.  This implies that
Lk$(K_i,l_j^t)=-$Lk$(K'_i,l_j^t)$ for all $1\leq j \leq \alpha $
and $1 \leq t \leq p$.\\
Consequently the linking matrix of $L$ is of the form (here
coefficients are considered modulo $p$).

$$\left (\begin{array}{cccccccc}
0&0&\dots&0& t_2^1&t_2^2&\dots&t_{\alpha}^p\\
0&0&\dots&0& -t_2^1&-t_2^2&\dots&-t_{\alpha}^p\\
0&0&\dots&0&.&.&\dots&.\\
0&0&\dots&0&.&.&\dots&.\\
t_2^1&-t_2^2&\dots&.&.&.&\dots&.\\
t_2^2&-t_2^2&\dots&.&.&.&\dots&.\\
.&.&\dots&.&.&.&\dots&.\\
.&.&\dots&.&.&.&\dots&.\\
.&.&\dots&.&.&.&\dots&.\\
t_{\alpha}^p&-t_{\alpha}^p&\dots&.&.&.&\dots&.
\end{array} \right )
$$

It is well known that the lines of this linking matrix  sum to
zero. In addition, the first line and the second line of our
matrix are dependent. Hence, all the cofactors of the matrix are
null. This concludes that $a_0$ is null modulo
$p$. \fin\\
\\

\textbf{Proposition 4.1.} {\sl If all the links of type $m$ have
$a_{2i-2}$  null modulo $p$, where $2i < p-1$. Then all links of
type $m-1$ have $a_{2i}$ null
modulo $p$.}\\
\emph{Proof.} Assume that for all links $L''$ of type $m$, we have
$a_{2i-2}(L'')\equiv 0$   modulo $p$, where  $2i < p-1$. Let $L$
be a link of type $m-1$. By definition we have $L=K_1 \cup
K_1'\cup \dots \cup K_{m-1}\cup K_{m-1}' \cup L'$, where $L'$ is a
strongly periodic OS link.
 We start by proving  the following lemma\\

\textbf{Lemma 4.4.} {\sl If $L_{p+}$ and $L_{p-}$ are two
$p$-periodic links of type $m-1$  such that their quotients differ
only by a self-crossing change in $\overline {L'}$. Then
$a_{2i}(L_{p+}) \equiv
a_{2i}(L_{p-}) $  modulo $p$.}\\

\emph{ Proof.} If we change a self-crossing   in the quotient link
$\overline {L'}$, then according to Lemma 2.1  the link $L_{p0}$
is either of type $m$ having $(\sharp L -p+2)$ components or a
periodic link having  $(\sharp L +p)$ components. In the second
case Remark 4.1 implies that $a_{2i}(L_{p+}) \equiv a_{2i}(L_{p-})
$ modulo $p$. However in the first case, we get
 $a_{2i} (L_{p+})-a_{2i}(L_{p-})\equiv a_{2i-2}(L_{p0})$ modulo
$p$. Since $L_{p0}$ is of type $m$ then $a_{2i-2}(L_{p0})$ is null
modulo $p$. This ends the proof of  Lemma 4.4. \fin \\
\\

\textbf{Lemma 4.5.}  {\sl Let $i$ be an integer such that
$2i<p-1$. If $a_{2i}(L) \equiv 0$ modulo $p$ for links of type
$m-1$ with
algebraically split orbits. Then, $a_{2i}(L) \equiv 0$ for  links of type $m-1$.}\\
\emph{Proof.} The idea of the proof is to use self-crossing
changes to split the orbits of our link. The algorithm can be
described as
 follows. We start by changing the crossings of $l_1^1$ and
 $l_1^2$, until having their linking number zero. Of course, at the same
time we get $lk(\varphi_p^i(l_1^1),\varphi_p^i(l_1^2))=0,$ for all
$1\leq i \leq p-1 $. In the second step we do the same for $l_1^1$
and $l_1^3$ and so on. We apply the same procedure for each orbit
of components until getting all orbits algebraically
split.  \fin \\
\\

Return now to the proof of Proposition 4.1. According to Lemma
4.5, it will be enough to prove the proposition in the case where
$L$ has only algebraically split orbits. The proof is done by
induction on $\alpha$. Remember that the number of components of
$L$ is $2(m-1)+p\alpha$. If
$\alpha=0$ then $a_{2i}$ is zero modulo $p$, due to Lemma 4.2.\\
Now assume that  the result is true for links of type $m-1$ with
less than $2(m-1)+p\alpha$ components. Assume that $L$ is a link
of type $m-1$ with orbits algebraically split such that $L$
has $2(m-1)+p\alpha$ components.\\
Since the orbits of $L$ are algebraically split we can use mixed
equivariant crossing changes to transform our link into an
algebraically split link. According to Lemma 2.1, there are 4
types of mixed equivariant crossing changes:\\
In the first and the second case of Lemma 2.1, Remark 4.1 implies
that $a_{2i}(L_{p+})-a_{2i}(L_{p-}) \equiv 0$ modulo $p$.\\
In the third and the fourth case, Lemma 4.1 implies that
$a_{2i}(L_{p+})-a_{2i}(L_{p-}) \equiv a_{2i}(L_{p0})$ modulo $p$.
Since the link $L_{p0}$ is of type $m-1$ having less components
than the original link $L$, then we can use the induction
 hypothesis to conclude that $a_{2i}(L_{p+})-a_{2i}(L_{p-}) \equiv 0$ modulo $p$.\\
In conclusion, mixed equivariant crossing changes do not affect
the coefficient $a_{2i}$ modulo $p$. Using these  changes we
should be able  to transform our link into an algebraically split
link having $a_{2i}$ null modulo
$p$. This completes the proof of Proposition 4.1. \fin\\
\\

The rest of  the proof of the Main Theorem is  a matter of
induction. According to Lemma 4.3, a link of type $(p-1)/2$ has
$a_{0}$ null modulo $p$. Starting from this fact, an easy
induction using Proposition 4.1 shows that a link of type zero
(which is a strongly periodic OS link) has $a_{2i}$ null modulo
$p$ for all
$i$ such that $2i< p-1$. This completes the proof of the Main Theorem. \fin\\


\begin{thebibliography}{99}
\bibitem{BM} H. {\sc Bass} and J. W. {\sc Morgan}. {\em The Smith
conjecture}. Pure and App. Math. 112, New York Academic Press
(1994).

\bibitem{Bu} G. {\sc Burde}. {\em \"{U}ber periodische knoten}. Arkiv der
math. (Basel) 30, (1998), pp. 487-492.


\bibitem{Ch1} N. {\sc Chbili}. {\em Les invariants $\theta_{p}$ des
3-vari\'et\'es p\'eriodiques}. Annales de l'institut Fourier,
Fascicule 4, (2001), pp.  1135-1150.

\bibitem{Ch2} N. {\sc Chbili}. {\em Quantum invariants and finite group actions on
 3-manifolds}, Topology
 Appl. vol 136/1-3, (2004), pp. 219-231.
 \bibitem{Ch3} N. {\sc Chbili}. {\em The Casson-Walker-Lescop invariant
  of periodic three-manifolds}, Math. Proc. Cambridge Phil. To appear
  Vol 140, 1, (2006).


\bibitem{CL} Q. {\sc Chen} and {\sc  T. Le}. {\em Quantum
invariants of periodic links and periodic 3-manifolds}. Fund.
Math. 184 (2004), 55-71.

\bibitem{DL} {\sc J. F. Davis,} and {\sc C. Livingston}. {\em
Alexander polynomials of periodic knots}. Topology, 30, (1991),
pp. 551-564.


\bibitem{BGP}  P. {\sc Gilmer}, J. {\sc Kania-Bartoszynska},
and J. {\sc Przytycki}. {\em 3-Manifold invariants and periodicity
of homology spheres}. Algebraic and Geometric Topology, Volume
2,(2002), pp. 825-842.
\bibitem{Ha} R. {\sc Hartley}. {\em The Conway potential function for links.} Comment. Math.
Helv. 58 (1983), no. 3, 365--378.

\bibitem{Hi} J. {\sc Hillman}. {\em New proofs of two theorems on periodic
knots}.  Archiv. Math.  37,  (1981), pp. 457-461.

\bibitem{Hoso} F. {\sc Hosokawa}. {\em On $\nabla $-polynomials of links.} Osaka
Math. J. 10 (1958) 273--282.

\bibitem{Ho} J.{\sc Hoste}. {\em The first
coefficient of the Conway polynomial}. Proc. Amer. Math. Soc., 95,
(1985), pp. 299-302.

\bibitem{Les}  C. {\sc Lescop}. {\em Global surgery formula for the Casson-Walker invariant}, Annals of
 Mathematics Studies, Princeton Univ.
 Press (1996).

\bibitem{Le} J. {\sc
Levine}. {\em The Conway polynomial of an algebraically split
link}. Proceeding of knots 96, edited by S. Suzuki, World
Scientific Publishing Co, (1997), pp. 23-29.

\bibitem{Mi} Y. {\sc Miyazawa}. {\em Conway polynomials of periodic links}. Osaka J. math.  31, pp. 147-163, 1994.
\bibitem{Mu} K. {\sc Murasugi}. {\em On periodic knots}. Comment.
Math. Helv. 46, (1971), pp. 162-174.


\bibitem{Pr} J. H. {\sc  Przytycki}. {\em On Murasugi's and Traczyk's criteria for periodic
links}. Math. Ann., 283, (1989), pp. 465-478.

\bibitem{PS} J.  {\sc Przytycki }
and M. {\sc Sokolov}. {\em Surgeries on periodic links and
homology of periodic 3-manifolds.} Math. Proc. Cambridge Phil.
Soc. Vol. 131(2), (2001), pp 295--307.

\bibitem{Sa} M. {\sc Sakuma}. {\em On the polynomials of  periodic links}.
 Math. Ann.  257, (1981), pp. 487-494.
\bibitem{Sa1} M. {\sc Sakuma}. {\em Surgery description of orientation-preserving periodic maps on
compact orientable 3-manifolds.} Rend. Istit. Mat. Univ. Trieste
32 (2001), suppl. 1, 375--396 (2002).

\end{thebibliography}
\end{document}